\documentclass{TEMA}

\usepackage[brazil]{babel}      % para texto em Português

\usepackage[latin1]{inputenc}   % para acentuação em Português
\usepackage[pdftex]{graphicx}
\usepackage{subfigure}
\usepackage{graphicx}
\usepackage{psfrag}
\usepackage{epsfig}
\usepackage{float}

\begin{document}

%********************************************************
\title
    {Formula\c c\~oes Semi-Discretas para a Equa\c c\~ao 1D de Burgers}

\author
    {C.A. LADEIA%,
    \thanks{cibele\_mat\_uel@yahoo.com.br; Doutoranda da Universidade Federal do Rio Grande do Sul
      }\,,
     N.M.L. ROMEIRO%,
    \thanks{nromeiro@uel.br;}\,, P.L. NATTI%,
    \thanks{plnatti@uel.br;}\,, E.R. CIRILO%,
    \thanks{ercirilo@uel.br;}\,,
    Departamento de Matemática, UEL, 86051-990 Londrina, PR, Brasil
     %\\ \\
     % V.G. FERREIRA%,
     %\thanks{pvgf@icmc.usp.br;}\,,
     %Departamento de Matemática Aplicada e Estatística, ICMSC, 13560-970 São Carlos, SP, Brasil
     \\ \\}
\criartitulo

\runningheads {Ladeia, Romeiro, Natti e Cirilo}{Formula\c c\~oes Semi-discretas para a Equa\c c\~ao 1D de Burgers}

\begin{abstract}
{\bf Resumo}. Neste trabalho fizemos comparações entre formulações semi-discretas para a obtenção de soluções numéricas para a equação 1D de Burgers. As formulações consistem em
discretizar o domínio temporal via métodos implícitos
multi-estágios de segunda e  quarta ordem: aproximantes de Padé $R_{11}$ e $R_{22}$; e o domínio espacial via métodos de elementos finitos: mínimos quadrados (MEFMQ), Galerkin (MEFG) e \textit{Streamline-Upwind} Petrov-Galerkin (SUPG). Conhecendo as soluções analíticas da equação 1D de Burgues, para diferentes condições iniciais e de fronteira, foram realizadas  análises dos erros numéricos  a partir das normas $L_{2}$ e $L_{\infty}$. Verificamos que o método com o aproximante de Padé $R_{22}$ adicionado as formulações MEFMQ, MEFG e SUPG, aumentou a região de convergência das soluções numéricas e apresentou maior precisão quando comparado as soluções obtidas por meio do aproximante de Padé $R_{11}$. Constatamos que o método $R_{22}$ amenizou as oscilações das soluções numéricas associadas as formulações  MEFG e SUPG.

{\bf Palavras-chave}. Equação de Burgers, Aproximantes de Padé, Métodos implícitos multi-estágios,
Métodos de elementos finitos.
\end{abstract}

%********************************************************
\newsec{Introdução}
Com a evolução da mecânica computacional  ocorreu a intensificação de pesquisas na resolução numérica de equações diferenciais com impactos positivos para a sociedade \cite{Oden}. Devido à vasta aplicabilidade destas pesquisas em problemas que envolvem processos convectivos e difusivos \cite{DoneaHuerta,Valdemir1,Valdemir2,Jain,RomeiroActa,RomeiroLandau}, interessa-nos as formulações numéricas que possam ser aplicadas a estes problemas, em especial à equação de Burgers \cite{David,Dogan,DoneaHuerta, Donea,Kakuda,Burgers4,Tian,venutilli}.

Vários autores apresentaram soluções numéricas para a equação 1D de Burgers, usando métodos de elementos finitos  \cite{Dogan,DoneaHuerta,Donea,Jain,Jiang,Burgers4,Burgers6,Tian,venutilli}, de elementos de contorno \cite{Kakuda}, de diferenças finitas \cite{Donea,Valdemir2,Jiang}, assim como alternativas para resolver o termo não linear convectivo como esquemas  \textit{upwind} \cite{Valdemir1,Valdemir2} e técnicas de linearização \cite{Dogan,Jain,Burgers4}.
Neste contexto faremos comparações de formulações semi-discretas para a equação 1D de Burgers, onde utilizaremos formulações semi-discretas para a discretização temporal e espacial e realizaremos uma linearização no termo convectivo. Esta linearização altera o tamanho do elemento em cada etapa  usando a informação a partir do passo anterior \cite{Dogan,Jain,Burgers4} transformando a equação de Burgers em um problema linear local.

As formulações semi-discretas consistem em discretizar o domínio temporal utilizando métodos implícitos multi-estágios e o domínio espacial via métodos de ele\-mentos finitos \cite{Donea,Tian,venutilli}. Em particular, para resolver a equação 1D de Bur\-gers, consideramos os métodos implícitos multi-estágios de segunda ordem, $R_{11}$, e de quarta ordem, $R_{22}$, e três formulações de métodos de elementos finitos: mínimos quadrados (MEFMQ), Galerkin (MEFG) e {\it Streamline-Upwind} Petrov-Galerkin (SUPG).

 Como o objetivo deste trabalho é, conhecido a solução analítica da equação 1D de Burgers, determinar quais formulações semi-discretas  apresentam soluções numéricas com  precisão em um intervalo maior de tempo, aumentando desta forma a região de convergência, apresentamos resultados comparativos entre as soluções e simulações numéricas para diferentes tempos. Apresentamos também tabelas com os cálculos dos logaritmos dos erros entre as formulações utilizando as normas $L_{2}$ e $L_{\infty}$.

 Desta forma, este artigo encontra-se estruturado como segue: a seção 2 apresenta a equação 1D de Burgues na forma geral e também na forma padrão de operador; na seção 3 é apresentada a discretização espacial obtida por meio dos aproximantes de Padé; na seção 4 são apresentadas as formulações de elementos finitos, utilizadas na discretização espacial da equação 1D de Burgues; na seção 5 é introduzida a técnica de linearização do termo convectivo transformando a equação de Burgers em um problema linear local; a seção 6 mostra os resultados numéricos obtidos no trabalho e finalmente as conclusões são apresentadas na seção 7.

\section{Equação 1D de Burgers}
Seja a equação 1D de Burgers definida em um domínio $\Omega=[a, b] \subset R$  limitado e aberto com fronteira $\Gamma = \partial\Omega$,  satisfazendo
\begin{eqnarray}
\label{burgers}
&& u_{t}(x,t)+ u(x,t)u_{x}(x,t) -  \epsilon u_{xx}(x,t) = f(x,t) \ \ \ \hbox{em} \  \Omega,\\
\label{condicao}
&& u(x,0)=u_{0}(x) \ \ \ \ \ \ \ \ \ \ \ \ \ \forall \ x\in \Omega,\\
\label{direchlet2}
&& u(a,t)=c_{0}=u(b,t)\ \ \ \ \ \ \ \ \ \ \  \hbox{em} \ \Gamma,
\end{eqnarray}
onde $\epsilon = \frac{1}{Re}$, $Re$ é o número adimensional de Reynolds; $u(x,t)$ é a componente da velocidade do fluido na direção do eixo $x$ e $f(x,t)$ representa o termo fonte. A equação  $(\ref{condicao})$ é uma condição inicial com $u_{0}$ sendo uma função dada e $(\ref{direchlet2})$ é uma condição de fronteira do tipo Dirichlet, sendo $c_{0}$ uma constante.

Reescrevendo ($\ref{burgers}$) na forma  padrão de operador, temos
\begin{eqnarray}
\label{burgers2}
u_{t} + \mathcal{L}(u)  = f,
\end{eqnarray}
onde $\mathcal{L}(u)= uu_{x} - \epsilon u_{xx}$   representa o  operador espacial e  descreve a soma dos operadores não linear e linear, convectivo e difusivo, respectivamente.
%$\mathcal{L}=\mathcal{L}_{conv}+\mathcal{L}_{dif}$
%\begin{eqnarray}
%\label{burgers3}
%\mathcal{L}(u)= uu_{x} - \epsilon u_{xx}
%\end{eqnarray}
%deverão ser encaminhados em CD(s) (com uma cópia impressa em
%papel A4) ao endereço: \\ %[1.3ex]
%\hspace*{1.3cm}TEMA/Fluxo Contínuo \\%
%\hspace*{1.3cm}SBMAC \\ %
%\hspace*{1.3cm}a/c Prof. Dr. Edson Wendland \\
%\hspace*{1.3cm}Caixa Postal 668 \\
%\hspace*{1.3cm}13560-970   São Carlos, SP.\\[1.3ex]

%%%%%%%%%%%%%%%%%%%%%%%%%%%%%%%%%%%%%%%%%%%%%%%%%%%%%%%%%%%%%%%%%%%%%%
\section{Discretização temporal}
Muitas técnicas numéricas para a discretização temporal são utilizadas para resolver equações diferencias parciais \cite{Donea,Valdemir1,John}. Consideramos a técnica passo de tempo obtida por meio dos aproximantes de Padé, cuja aproximação para o operador de evolução é dada por
\begin{eqnarray}
\label{operadortempo}
\mathcal{E}(\Delta t): u(t^{n})\rightarrow u(t^{n+1}),
\end{eqnarray}
que permite transportar a solução numérica em um determinado tempo $t^{n}=n\Delta t$ para o próxi\-mo tempo $t^{n+1}=t^{n}+ \Delta t$, sendo  $\Delta t$ um passo de tempo. Nesta situação, a evolução de $\mathcal{E}(\Delta t)$,  obtida  a partir do desenvolvimento da série de Taylor, é dada por meio do operador exponencial, ou seja
\begin{eqnarray}
\label{exptempo}
u(t^{n+1})&=& \left(1 + \Delta t\frac{\partial}{\partial t}   + \frac{1}{2!}\Delta t^{2}\frac{\partial^{2}}{\partial t^{2}} + \frac{1}{3!}\Delta t^{3}\frac{\partial^{3}}{\partial t^{3}}+  ... + \frac{1}{n!}\Delta t^{n}\frac{\partial^{n}}{\partial t^{n}}+ \ldots\right )u(t^{n})\nonumber \nonumber\\
&=&\exp\left(\Delta t\frac{\partial}{\partial t}\right)u(t^{n}).
\end{eqnarray}
 Assim, a técnica de passo de tempo de várias ordens pode ser obtida usando os aproximantes de Padé para o operador exponencial $e^{h}$, sendo $h=\Delta t\frac{\partial}{\partial t}$ \cite{David,Donea}. Os primeiros aproxi\-mantes de Padé para $e^{h}$ encontram-se apresentados na Tabela 1.
\begin{table}[h]
\caption{Aproximantes de Padé da função exponencial $e^{h}$} \label{tabela01}
\begin{center}
\begin{tabular}{c| c c c c}
 \hline
 $R_{LM}$ & $M$=0 &  $M$=1 &  $M$=2 &  $M$=3 \\
  \hline
  $L$=0 & 1 & $\frac{1}{1-h}$& $\frac{2}{2-2h+h^{2}}$& $\frac{6}{6-6h+3h^{2}-h^{3}}$ \\
  \hline
  $L$=1 & $\frac{1+h}{1}$& $\frac{2+h}{2-h}$&$ \frac{6+2h}{6-4h+h^{2}}$& $\frac{24+6h}{24-18h+6h^{2}-h^{3}}$\\
  \hline
   $L$=2 & $\frac{2+2h+h^{2}}{2}$& $\frac{6+4h+h^{2}}{6-2h}$& $\frac{12+6h+h^{2}}{12-6h+h^{2}}$&$\frac{60+24h+3h^{2}}{60-36h-9h^{2}-h^{3}}$ \\
  \hline
   $L$=3 & $\frac{6+6h+3h^{2}+h^{3}}{6}$&$\frac{24+18h+6h^{2}+h^{3}}{24-6h}$ & $\frac{60+36h+9h^{2}+h^{3}}{60-24h+3h^{2}}$ &$\frac{120+60h+12h^{2}+h^{3}}{120-60h+12h^{2}-h^{3}}$ \\
  \hline
  \end{tabular}
  \end{center}
\end{table}

Em \cite{Hairer, Lambert} os autores mostram que a aproximação de Padé $R_{LM}$  é incondicionalmente estável se satisfaz a condição $M-2\leq L\leq M$. Assim, as técnicas de multi-estágios implícitas empregadas
na parte da discretização temporal, utilizadas neste trabalho, são incondicionalmente estáveis com erro de truncamento  de ordem $\mathcal{O}(\Delta t^{2n})$ \cite{David}. Logo, podemos reescrever o método implícito na seguinte forma de operador
\begin{eqnarray}
\label{tempo}
\frac{\Delta u}{\Delta t}- \mathbf{W}\triangle u_{t} = \textmd{\textbf{w}}u_{t}^{n},
\end{eqnarray}
onde o vetor $\Delta u$  tem $n$ componentes (ou $n_{\textit{estágios}}$) \cite{Donea} e $\Delta u_{t}$ é a derivada parcial de $\Delta u$ com respeito ao tempo. Substituindo ($\ref{burgers2}$)  em ($\ref{tempo}$) e considerando que o operador $\mathcal{L}$ é linear com coeficientes constantes, temos que
\begin{eqnarray}
\label{tempocdr1}
\frac{\Delta u}{\Delta t}+ \mathbf{W} \mathcal{L} (\Delta u) = \textmd{\textbf{w}}[f^{n}-\mathcal{L}(u^{n})] + \mathbf{W}\Delta f
\end{eqnarray}
com $\Delta u$, $\mathbf{W}$, $\Delta f$, e $\textmd{\textbf{w}}$ dependentes do método escolhido. Para os métodos implícitos $R_{11}$ e $R_{22}$, os valores de $\Delta u$, $\mathbf{W}$, $\Delta f$, e $\textmd{\textbf{w}}$ encontram-se definidos nas seções 3.1 e 3.2. Devido à limitação de espaço, definimos apenas as formulações compactas destes métodos, outros detalhes podem ser obtidos em  \cite{Donea,Huertabernadino}.\\

\subsection{Método de segunda ordem}
A formulação compacta  $R_{11}$, Crank-Nicolson  \cite{Donea,Huertabernadino}, é dada por
\begin{eqnarray}
\label{cn}
&&\Delta u = u^{n+1} - u^{n}, \ \ \ \ \Delta f = f^{n+1} - f^{n}\nonumber\\
&&\mathbf{W} = 1/2, \ \ \ \ \ \ \ \ \ \ \ \ \ \ \   \textmd{\textbf{w}}=1.
\end{eqnarray}
%Neste caso a dimensão  é $1$ e consequentemente os vetores e matrizes são %escalares.\\
\subsection{Método de quarta ordem}
A formulação compacta $R_{22}$ \cite{Donea,Huertabernadino} é
\begin{eqnarray}
\label{r22}
&&\Delta u =
\left\{
\begin{array}{cc}
 u^{n+1/2} - u^{n} \\
u^{n+1} - u^{n+1/2}
\end{array}
\right\},
 \ \ \ \
\Delta f =
\left\{
\begin{array}{cc}
 f^{n+1/2} - f^{n} \\
f^{n+1} - f^{n+1/2}
\end{array}
\right\},\nonumber\\
 &&\mathbf{W} =\frac{1}{24}  \left[
\begin{array}{cc}
7 & -1 \\
13 & 5
\end{array}
\right],  \ \ \ \ \ \ \ \ \ \ \ \ \
\textmd{\textbf{w}}= \frac{1}{2} \left\{
\begin{array}{cc}
1 \\
1
\end{array}
\right\}.
\end{eqnarray}

\section{Discretização  espacial}

Neste trabalho utilizamos o método de elementos finitos para a discretização espacial. Este método introduz funções bases $\{\varphi_{0},\ldots,\varphi_{m}\}$ que geram o subes\-paço onde está sendo procurada a solução exata, com suporte locali\-zado  nos pontos nodais dos elementos \cite{Jiang}. Para definir as funções bases realizamos uma discretização no intervalo $[a,b]$, dividindo-o em $m$ sub-intervalos (ou elementos) $e_{j}=(x_{j-1},x_{j})$, $j=0,1,\ldots,m$, de comprimento $h_{j}= x_{j}- x_{j-1}.$

Seja $u_{h}$ uma função teste linear em cada elemento $e_{j}$ e contínua sobre $[a,b]$ satisfazendo as condições de fronteira $u_{h}(a)=c_{0}=u_{h}(b)$. Escolhendo os valores $u_{0}, u_{1}, u_{2},\ldots,u_{m}$ nos nós $x_{j}$,  como parâmetros para descrever $u_{h}(x)$ e expressando-os sobre cada elemento $e_{j}$, obtemos $u_{h}(x)=\psi_{1}^{(j)}(x)u_{j-1}+\psi_{2}^{(j)}(x)u_{j}$ \ com \ $x\in e_{j}$,
%begin{eqnarray}
%label{mef}
%_{h}(x)=\psi_{1}^{(j)}(x)u_{j-1}+\psi_{2}^{(j)}(x)u_{j} \ \ x\in e_{j},
%end{eqnarray}
onde
\begin{eqnarray}
\label{mef1}
\psi^{(j)}(x)= \left(\begin{array}{cc} \psi_{1}^{(j)}\\ \psi_{2}^{(j)}
\end{array}\right)= \left(\begin{array}{cc} (x_{j}-x)/h_{j} \\ (x-x_{j-1})/h_{j} \end{array}\right).
\end{eqnarray}
Podemos escrever $u_{h}(x)$, sobre todo o domínio $[a,b]$, como sendo
\begin{eqnarray}
\label{mef2}
u_{h}(x)= \varphi_{0}(x)u_{0}+\varphi_{1}(x)u_{1}+ \ldots + \varphi_{m}(x)u_{m},
\end{eqnarray}
onde
\begin{eqnarray}
\label{52}
&&\varphi_{j}(x)=\left \{ \begin{array}{cc} \psi_{2}^{(j)}(x), \ \ \ \ \ \ \ \ \ \ \ \ \ \ \ \ \ \ \ \ \ \ \    x\in e_{j}; \\
                         \psi_{1}^{(j+1)}(x),\ \ \  x\in e_{j+1} \ \ \  1 \leq j \leq m-1 \\
                       \ \   0,\ \ \ \ \ \ \ \ \ \ \ \ \ \ \ \ \ \ \ \      \hbox{caso contrário};   \end{array}\right.\nonumber\\
\label{53}
&&\varphi_{0}(x)=\left \{ \begin{array}{cc} \psi_{1}^{(1)}(x), \ \ \ \ \ \ \ \ \ \ \ \ \ \ \ \ \ \ \ \ \ \ \    x\in e_{1}; \\
                         \ \   0,\ \ \ \ \ \ \ \ \ \ \ \ \ \ \ \ \ \ \ \      \hbox{caso contrário};   \end{array}\right.\nonumber\\
\label{54}
&&\varphi_{n}(x)=\left \{ \begin{array}{cc} \psi_{2}^{(m)}(x), \ \ \ \ \ \ \ \ \ \ \ \ \ \ \ \ \ \ \ \ \ \ \    x\in e_{m}; \\
                         \ \   0,\ \ \ \ \ \ \ \ \ \ \ \ \ \ \ \ \ \ \ \      \hbox{caso contrário};   \end{array}\right.\nonumber
\end{eqnarray}
são as funções bases e $\varphi_{j}(x)$ é uma função linear seccionalmente contínua, com valor um para o nó $x_{j}$ e nulo para os outros nós. Portanto, $u_{h}(x_{j})=u_{j}$. Especificamente, se temos $u_{0}=c_{0}$ então ($\ref{mef2}$) satisfaz as condições de fronteira ($\ref{direchlet2}$) e para os parâmetros $u_{1},u_{2},\ldots,u_{m}$ valores arbitrários. Todas as funções $u_{h}(x)$ constituem o espaço das funções teste $V_{h}$. Considere o espaço de Hilbert $V=H^{1}_{0}(\Omega)$, onde $V$ é o conjunto das soluções tentativas. Assim, seja $V_{h}$  o subespaço de dimensão finita do espaço de dimensão infinita $V,$  formado por funções lineares seccionais geradoras de um conjunto de $m$ elementos de $V$ denotado por $V_{h}=[\varphi_{0},\ldots,\varphi_{m}]$. As funções bases $\varphi_{j}$ são obtidas a partir do método de elementos finitos, considerando a partição $a=x_{0}<x_{1}<x_{2} \ldots <x_{m-1}<x_{m}=b.$

\subsection{Método de elementos finitos via  mímimos quadrados}
O método de elementos finitos via  mímimos quadrados (MEFMQ) consiste em aproximar os termos  da parte espacial de ($\ref{burgers}$) por meio de uma  formulação variacional  \cite{Jiang}, obtida usando o método de mínimos quadrados que minimiza o quadrado da integral do resíduo.

Seja $V$ o espaço de Hilbert e defina o funcional
\begin{eqnarray}
\mathcal{F}: V &\rightarrow& R \nonumber\\
u^{j+1}&\rightarrow& \mathcal{F}(u^{j+1}),
\end{eqnarray}
onde $\mathcal{F}(u)=\int_\Omega [u(x)u_x(x)-\epsilon u_{xx}(x)-f(x)]^2 dx$ para todo $x\in \Omega$ sobre todos os $u\in V$. Minizando  $\mathcal{F}$ com respeito a $u^{j+1}$ para $j=0,1,2,\ldots,m$, e usando a derivada de Gâteaux \cite{Daryoush e Encyeh}, podemos resolver  o seguinte problema variacional: achar $u^{j+1} \in V$ tal que
\begin{eqnarray}
a_{M}(u^{j+1}; u_{x}^{j+1}, w)= F_{M}(u;w), \ \ \ \forall  w \in V,
\end{eqnarray}
onde definimos o funcional bilinear $a_{M}(u;\cdot,\cdot): V\times V \rightarrow R$ e o funcional linear $F_{M}(u;\cdot): V \rightarrow R$ por \cite{Cibele}
\begin{eqnarray}
&&a_{M}(u; u_{x},w) = \int_\Omega[\epsilon^{2}u_{xx}w_{xx}-\epsilon uu_{x}w_{xx}]dx, \ \ \ \forall w \in V; \\
&&F_{M}(u;w)= \int_\Omega [-\epsilon w_{xx}f + u_{x}wf + uw_{x}f]dx, \ \ \ \forall w \in V.
\end{eqnarray}
Obtido a formulação variacional, podemos resolver a parte espacial de ($\ref{burgers}$) utilizando MEFMQ, considerando o subespaço $V_{h}\subset V$, que consiste em: determinar uma solução aproximada $u^{j+1}_{h} \in V_{h}$ tal que
\begin{eqnarray}
a_{M}(u^{j+1}_{h};v_{h}^{j+1}, w_{h})= F_{M}(u_{h};w_{h}), \ \ \ \forall w_{h}
 \in V
\end{eqnarray}
onde $v = u_{x}$.

\subsection{Método elementos finitos via Galerkin}
Uma segunda aproximação para a parte espacial de ($\ref{burgers}$) é realizada através do método de Galerkin (MEFG) \cite{Brenner e Scott,Ciarlet}. Seja $V$ o espaço de Hilbert e dada a formulação fraca, podemos resolver o seguinte problema: determinar $u^{j+1} \in V$ tal que
\begin{eqnarray}
\label{galer}
a_{G}(u^{j+1}; u_{x}^{j+1}, w)= F_{G}(w), \ \ \ \forall w \in V,
\end{eqnarray}
onde definimos o funcional bilinear $a_{G}(u;\cdot,\cdot): V\times V \rightarrow R$ e o funcional linear $F_{G}(\cdot): V \rightarrow R$ por \cite{Cibele}:
\begin{eqnarray}
&&a_{G}(u; u_{x},w) = \int_\Omega [uu_{x}w-\epsilon u_{x}w_{xx}]dx, \ \ \ \forall w \in V; \\
&&F_{G}(w)= \int_\Omega wf dx, \ \ \ \forall w \in V.
\end{eqnarray}
Para resolver a parte espacial de  ($\ref{burgers}$), utilizando MEFG, consideramos $V_{h}\subset V$, no qual o problema agora consiste em determinar uma solução aproximada $u^{j+1}_{h} \in V_{h}$ tal que
\begin{eqnarray}
\label{1galerkin}
a_{G}(u^{j+1}_{h}; v_{h}^{j+1}, w_{h})= F_{G}(w_{h}), \ \ \ \forall w_{h}\in V_{h}.
\end{eqnarray}

\subsection{Método estabilizado \textit{Streamline-Upwind} Petrov-Galerkin}
O método estabilizado \textit{Streamline-Upwind} Petrov-Galerkin (SUPG) contorna as limitações do método de Galerkin \cite{Brooks}. O SUPG é uma combinação da formulação de Galerkin com termos baseados no resíduo, em nível de elementos. Estes termos são balanceados por parâmetros de estabilização, resultando em formulações variacionais consistentes com as propriedades de estabilização, superiores às da aproximação de Galerkin \cite{Brooks,DoneaHuerta}.

Assim, o método estabilizado SUPG para aproximar a parte espacial  de ($\ref{burgers}$), consiste em determinar $u_{h}\in V_{h}$ tal que
\begin{eqnarray}
 \label{supg}
 a_{G}(u_{h}; v_{h},w_{h})+ E_{\textmd\small{SUPG}}(u_{h};v_{h},w_{h})= F_{G}(w_{h}), \ \ \ \forall w_{h} \in V_{h},
 \end{eqnarray}
onde $\frac{\partial u_{h}}{\partial x}=v_h$ e $E_{\textmd\small{SUPG}}(u_{h};v_{h},w_{h})$ indicam os termos de perturbação adicionados à formulação variacional padrão ($\ref{1galerkin}$). Estes termos são adicionados de forma a preservar a consistência do método para obter a estabilidade numérica, dada pela expressão
\begin{eqnarray}
 \label{supg2}
 E_{\textmd\small{SUPG}}(u_{h};v_{h},w_{h})&=& \sum_{e_{j}}\int_{e_{j}}u_{h}\frac{\partial w_{h}}{\partial x}\tau\left(u_{h} \frac{\partial u_{h}}{\partial x} - \epsilon\frac{\partial^{2} u_{h}}{\partial x^{2}} - f \right) d\Omega \nonumber\\
 &=&\sum_{e_{j}}(\mathcal{P}(u_{h};w_{h}),\tau \mathcal{R}(u_{h}))_{\Omega^{j}}
 \end{eqnarray}
sendo $\mathcal{P}(w)=u_{h}\frac{\partial w_{h}}{\partial x}$ a perturbação da função teste, enquanto o termo residual  $\mathcal{R}$ e o parâmetro $\tau$ são definidos, respectivamente por \cite{DoneaHuerta}:
\begin{eqnarray}
 \label{residualburgers}
 &&\mathcal{R}(u_{h})= u_{h}\frac{\partial u_{h}}{\partial x} -\epsilon\frac{\partial^{2} u_{h}}{\partial x^{2}} - f,\\
 \label{supgtau}
 &&\tau = \left(\left(\frac{2u}{h}\right)^{2} + 9\left(\frac{4\epsilon}{h^{2}}\right)^{2} \right)^{-1/2}.
 \end{eqnarray}

Para resolver a parte espacial do problema  ($\ref{burgers}$) utilizando SUPG, consideramos $V_{h}\subset V$ para $n=0,1,2,\ldots,N$, que consiste em determinar uma solução aproximada $u^{n+1}_{h} \in V_{h}$ tal que
\begin{eqnarray}
a_{G}(u^{j+1}_{h}; v^{j+1}_{h}, w_{h})+ E_{\textmd\small{SUPG}}(u^{j+1}_{h}; v^{j+1}_{h}, w_{h})= F_{G}(w_{h}), \ \ \ \forall w_{h}\in V_{h}.
\end{eqnarray}

\section{Linearização do termo convectivo}
 Várias técnicas para resolver o termo convectivo da equação de Burgers podem ser encontradas na literatura \cite{Dogan,DoneaHuerta, Valdemir1,Jain,Burgers4}. Neste trabalho, realizamos uma linearização no termo convectivo de ($\ref{burgers}$), que consiste em alterar o tamanho do elemento em cada etapa  utilizando a informação do passo anterior \cite{Dogan,Jain,Burgers4}, transformando a equação de Burgers em um problema linear local. Para isto, multiplicamos ambos os lados de ($\ref{burgers}$) por uma função teste $w \in V$ e integramos, resultando em
\begin{eqnarray}
\label{termoconv}
\int_\Omega (u_{t}+ uu_{x} - \epsilon u_{xx}- f)wdx=0.
\end{eqnarray}
Como a solução de $(\ref{burgers}$)-($\ref{condicao}$) é procurada  sobre o domínio $a\leq x\leq b$, com condições de fronteira em $x = a$ e em $x = b$, consideramos o subespaço de dimensão finita $V_{h}\subset V,$ onde as funções bases $\varphi_{j}$  são obtidas a partir do método elementos finitos, utilizando a partição de tamanho $h_{j}=x_{j}-x_{j-1}$, mapeada por uma coordenada local $\sigma$, onde $x_{j}=x_{j-1}+ \sigma h_{j}$, $0\leq\sigma\leq 1$  \cite{Dogan}. Logo, construímos uma função teste $u_{h}$ e escolhemos os valores $u_{0},u_{1}, u_{2}\ldots, u_{m}$ nos nós $x_{j}$. Podemos então reecrever a equação ($\ref{termoconv}$)  como
\begin{eqnarray}
\label{linearconvectivo2}
&\sum_{j=0}^{m}\int_\Omega \left(\frac{\partial u_{j}}{\partial t}+ \eta \frac{\partial \varphi_{j}(x) }{\partial x}u_{j} - \zeta \frac{\partial^{2} \varphi_{j}(x) }{\partial x^{2}}u_{j}- (\zeta - \eta) fu_{j} \right)\varphi_{i}(x)dx=0,
\end{eqnarray}
$\forall \ \varphi_{i}(x), \varphi_{j}(x) \in V_{h}$, onde $\eta =\frac{u_{0}}{h_{j}}$, $\zeta=\frac{\epsilon}{h^{2}_{j}}$ são localmente constantes sobre cada elemento e  $w_{h}=\varphi_{i}(x)$, para $i=0,1,2,\ldots,m,$ é definida como uma função teste.

Convém observar que neste trabalho a parte temporal de ($\ref{linearconvectivo2}$) é discretizada utilizando os aproximantes de Padé, $R_{11}$ e $R_{22}$,  e  a parte espacial via formulações dos métodos de elementos finitos: MEFMQ, MEFG e SUPG.

\section{Resultados Numéricos}
Todos os resultados apresentados na sequência resultam das formulações semi-discretas abordadas. Apresentamos  análises dos erros numéricos a partir das normas $L_{2}$ e $L_{\infty}$, comparando as soluções numéricas com as soluções analíticas dos exemplos avaliados.

\subsection{Exemplo 1}
A equação 1D de Burgers definida em (\ref{burgers}) com condição inicial
\begin{eqnarray}
\label{condburgers2}
u(x,0)=\frac{2\epsilon\pi\sin(\pi x)}{k + \cos(\pi x)}, \ \ \ k>1,
\end{eqnarray}
 e  condições de fronteira  $u(0,t)=0=u(1,t)$ tem como solução analítica \cite{Burgers4}
\begin{eqnarray}
\label{analburgers2}
u(x,t)=\frac{2\epsilon\pi\exp(\pi^{2}\epsilon t)\sin(\pi x)}{k + \exp(-\pi^{2}\epsilon t)\cos(\pi x)}, \ \ \ \hbox{com} \ k>1.
\end{eqnarray}

Considerando o domínio $0\leq x\leq 1$ e $k = 2$, apresentamos nas Figuras 1a-1f os resultados das formulações semi-discretas comparados com o resultado analítico, para os tempos $t=0.5$ e $t=1$, com  $\Delta t = 0.5$, $Re = 10^{5}$ e uma malha de 50 elementos lineares, o que equivale a $h = 0.02$.

\begin{figure}[!htb]
\centering
\epsfig{file=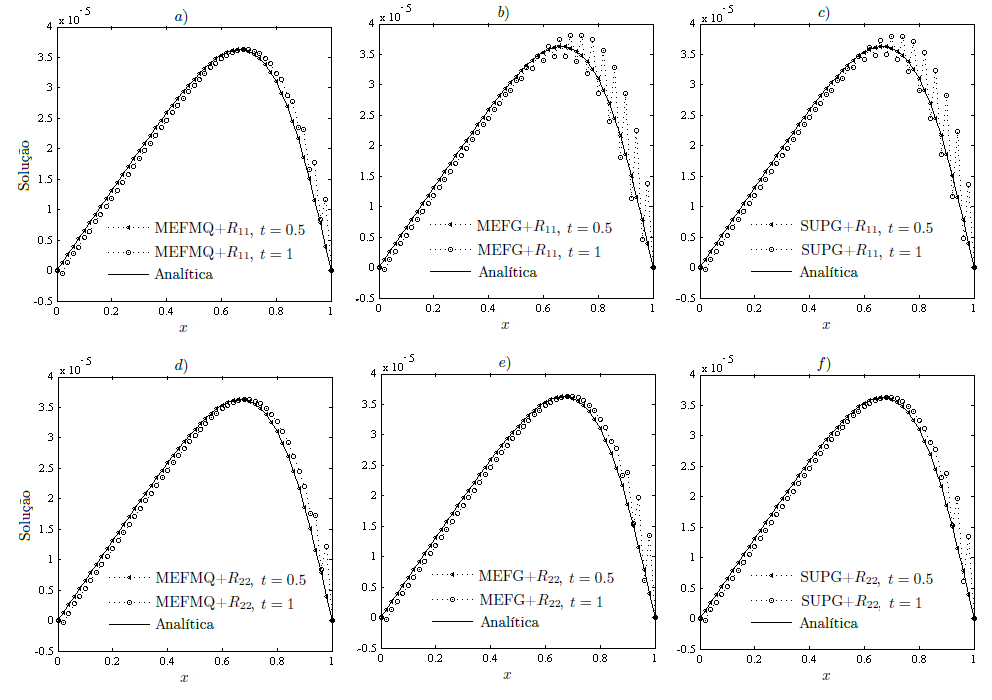,scale=0.49}
\caption{a-f) Solução analítica e soluções numéricas das formulações semi-discretas.} %
\label{fig1}
\end{figure}

Observamos que o aproximante de Padé $R_{11}$ adicionado as formulações MEFMQ, MEFG e SUPG, apresentou oscilações mais evidentes no tempo final $t=1$ conforme Figu\-ras 1a-1c, mas ressaltamos a ocorrência de menor intensidade para o caso \linebreak MEFMQ+$R_{11}$, Figura 1a. Ao incluirmos o aproximante de Padé $R_{22}$ nas formulações MEFMQ, MEFG e SUPG, verificamos  que as oscilações foram amenizadas, Figuras 1d-1f.

%É fato que o refinamento no tempo permite algum controle sobre as oscilações, porém a idéia deste exemplo é utilizar um $\Delta t$ grosseiro para evidenciar as oscilações.

Para uma avalição mais precisa das formulações semi-discretas, calculamos os logaritmos dos erros considerando  malhas com $h=1/50$ e $h=1/1000$ para os tempos $t=0.5$ e $t=1$ com  $\Delta t = 0.5$ e $Re =10^{5}$, tais resultados encontram-se na Tabela 2.

  \begin{table}[!htb]
\label{cdrnumerico1}
\centering
\caption{Cálculos dos logaritmos dos erros entre as formulações semi-discretas utilizando as normas $L_{2}$ e $L_{\infty}$.}
 \begin{tabular}{c c c c c c}
\multicolumn{6}{c}{ \ \ \ \ \ \ \ \ \ \ \ \ \ \ \ \ \ \ \ \ \ \ \ \ \ \ \ \ \ \ \ \ $h =1/50$\ \ \ \ \ \ \ \ \ \ \ \ \ \ \ \ \ \ \ \ \ \              $h =1/1000$  }\\
 \hline
 Formulações & E($\log_{10}$)       &$t=0.5$      & $t=1$               &$t=0.5$       & $t=1$   \\
 \hline
   MEFMQ+$R_{11}$ &$\|E\|_{2}$       &-1.4338   &-1.0964            &-1.4709    &-1.1762\\%\cline{2-4}
                 &$\|E\|_{\infty}$  &-1.0570   &-6.6807x$10^{-1}$  &-1.2759	   &-9.7797x$10^{-1}$ \\%\cline{3-5}
  \hline
  MEFMQ+$R_{22}$ &$\|E\|_{2}$       &-1.4363   &-1.0998           &-1.4649    &-1.1698 \\%\cline{2-4}
                 &$\|E\|_{\infty}$  &-1.0562   &-6.4226x$10^{-1}$  &-1.2566    &-9.5395x$10^{-1}$ \\
    \hline
   MEFG+$R_{11}$  &$\|E\|_{2}$       &-1.4117  &-8.6034x$10^{-1}$  &-1.4614    &-1.1665 \\%\cline{2-4}
                 &$\|E\|_{\infty}$  &-1.0112  &-5.1914x$10^{-1}$  &-1.1212	  &-8.1722x$10^{-1}$ \\
  \hline
  MEFG+$R_{22}$  &$\|E\|_{2}$       &-1.4289  &-1.0365            &-1.4612    &-1.1664 \\%\cline{2-4}
                 &$\|E\|_{\infty}$  &-1.0261  &-5.8131x$10^{-1}$  &-1.1178	  &-8.1754x$10^{-1}$ \\
  \hline
  SUPG+$R_{11}$  &$\|E\|_{2}$       &-1.4130   &-8.7854x$10^{-1}$ &-1.4614   &-1.1665\\%\cline{2-4}
                 &$\|E\|_{\infty}$  &-1.0128   &-5.2995x$10^{-1}$ &-1.1212	 &-8.1722x$10^{-1}$ \\
   \hline
  SUPG+$R_{22}$  &$\|E\|_{2}$       &-1.4286   &-1.0351           &-1.4615    &-1.1665 \\%\cline{2-4}
                 &$\|E\|_{\infty}$  &-1.0260   &-5.8146x$10^{-1}$ &-1.1514    &-8.3593x$10^{-1}$\\
     \hline
  \end{tabular}
\end{table}

A partir dos resultados apresentados na Tabela 2, confirmamos que o aproximante de Padé $R_{22}$ em conjunto com as formulações  MEFMQ,  MEFG e SUPG aumentou a região de convergência das soluções numéricas e apresentou maior precisão quando comparado as soluções obtidas por meio de aproximantes de Padé $R_{11}$, considerando $h =1/50$. Melhorando o refinamento na direção espacial, isto, é, tomando $h =1/1000$, observamos que a formulação MEFMQ+$R_{11}$  mostrou-se mais precisa entre todas as formulações, não havendo então a necessidade de um método oneroso para a discretizaçao no tempo.

\subsection{Exemplo 2}
Seja um problema de propagação uniforme de choque \cite{Burgers6}, para a equação 1D de Burgers (\ref{burgers}), com  condição inicial dada por
\begin{eqnarray}
\label{condburgers}
u(x,0)=\left\{
\begin{array}{rcl}
u_{1}=0.5,& \mbox{se} & x< 0\\
u_{2}=1.5, & \mbox{se} & x>0,
\end{array}
\right.
\end{eqnarray}
cuja solução analítica é dada por \cite{Burgers6}
\begin{eqnarray}
\label{analburgers}
u(x,t)= u_{1} + \frac{u_{2}-u_{1}}{1 + \exp\left[-Re\frac{u_{2}-u_{1}}{2}(x - \frac{u_{1}+u_{2}}{2}t)\right]}.
\end{eqnarray}

Considerando o domínio $-0.5\leq x\leq 0.5$  e as condições de fronteira satisfazendo $u(-0.5,0)=0.5$ e $u(0.5,0)=1.5$, apresentamos nas Figuras 2a-2f os resultados das formulações semi-discretas comparados com o resultado analítico, para os tempos $t=0.05$ e $t=0.1$, com $\Delta t = 3.3\times10^{-4}$, $Re =10^{4}$ e uma malha com 3000 elementos lineares.

\begin{figure}[!htb]
\centering
\epsfig{file=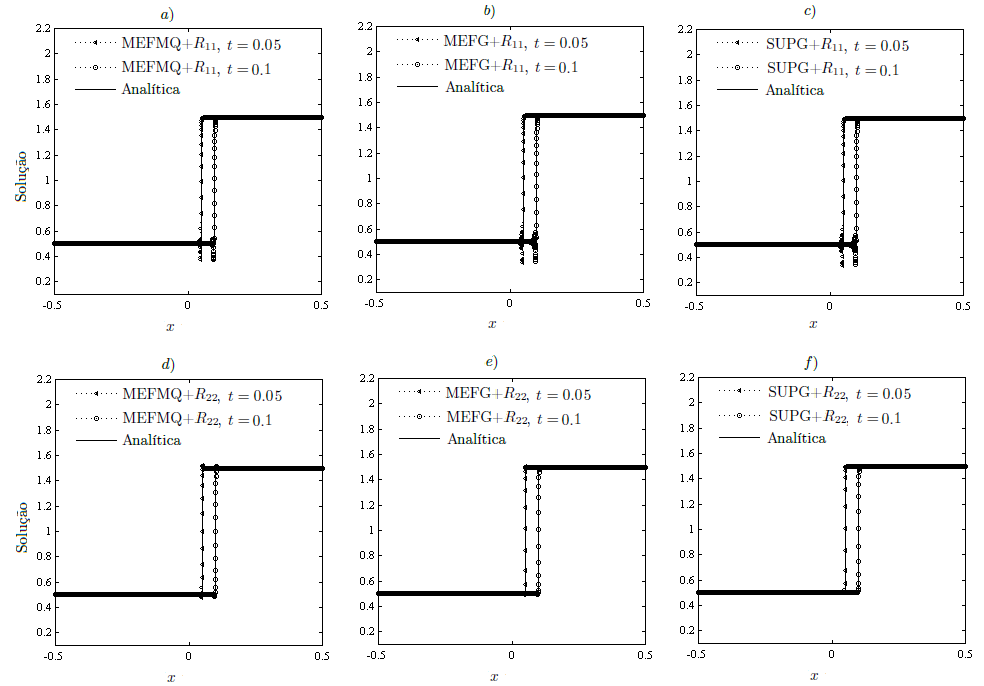,scale=0.49}
\caption{a-f) Soluções numéricas das formulações semi-discretas para  $h=1/3000$, $\Delta t = 3.3\times10^{-4}$,  $Re = 10^{4}$ e para os tempos $t=0.05$ e $t=0.1$.} %
\label{fig1}
\end{figure}

Podemos observar nas Figuras 2d-2f que ao adicionarmos o aproximante de Padé $R_{22}$ nas formulações
MEFMQ, MEFG e SUPG, obtivemos um aumento na região de convergência das soluções
numéricas, assim como uma maior precisão das soluções quando comparadas as obtidas por meio do aproximante de Padé $R_{11}$, vide Figuras 2a-2c.

\begin{table}[!htb]
\label{cdrnumerico1}
\centering
\caption{Cálculos
dos logaritmos dos erros entre as formulações semi-discretas utilizando as normas $L_{2}$ e $L_{\infty}$}
\tabcolsep=2pt
\begin{tabular}{c c c c c c}
\multicolumn{6}{c}{\ \ \ \ \ \ \ \ \ \ \ \ \ \ \ \ \ \ \ \ \ \ \ \ \ \ \ \ \ \  $h =1/3000$ \ \ \ \ \ \ \ \ \ \ \ \ \ \ \ \ \ \ \ \ \  $h =1/4000$ } \\
 \hline
 Formulações & E($\log_{10}$)   &$t=0.05$           & $t=0.1$            &$t=0.05$           & $t=0.1$ \\
\hline
  MEFMQ+$R_{11}$ &$\|E\|_{2}$   &-1.9468            &-1.8369             &-1.8382            &-1.7146 \\%\cline{3-5}
               &$\|E\|_{\infty}$&-6.5205x$10^{-1}$  &-5.8740x$10^{-1}$   &-6.4089x$10^{-1}$  &-5.9775x$10^{-1}$ \\
  \hline
  MEFMQ+$R_{22}$  &$\|E\|_{2}$  &-2.1093            &-1.9833             &-1.8465            &-1.7192  \\
               &$\|E\|_{\infty}$&-8.6062x$10^{-1}$  &-7.7387x$10^{-1}$   &-6.8386x$10^{-1}$  &-6.2629x$10^{-1}$ \\%\cline{3-5}
\hline
MEFG+$R_{11}$ &$\|E\|_{2}$      &-1.9063            &-1.8109             &-1.8397            &-1.7156 \\%\cline{3-5}
              &$\|E\|_{\infty}$ &-6.2961x$10^{-1}$  &-5.7230x$10^{-1}$   &-6.3991x$10^{-1}$  &-5.9758x$10^{-1}$\\
\hline
  MEFG+$R_{22}$  &$\|E\|_{2}$   &-2.2371            &-2.0699             &-1.8497            &-1.7206   \\%\cline{3-5}
               &$\|E\|_{\infty}$&-9.1682x$10^{-1}$  &-8.4416x$10^{-1}$   &-6.8569x$10^{-1}$  &-6.2722x$10^{-1}$\\
\hline
  SUPG+$R_{11}$  &$\|E\|_{2}$   &-1.9063            &-1.8109             &-1.8397            &-1.7156\\
              &$\|E\|_{\infty}$ &-6.2961x$10^{-1}$  &-5.7230x$10^{-1}$   &-6.3991x$10^{-1}$  &-5.9758x$10^{-1}$  \\
  \hline
  SUPG+$R_{22}$ &$\|E\|_{2}$    &-2.2473            &-2.0729             &-1.8496            &-1.7206 \\
             &$\|E\|_{\infty}$  &-9.4179x$10^{-1}$  &-8.5595x$10^{-1}$   &-6.8579x$10^{-1}$  &-6.2725x$10^{-1}$  \\
\hline
  \end{tabular}
\end{table}

Porém, confrontando  MEFMQ+$R_{22}$ com as formulações  MEFG+$R_{22}$ e \linebreak SUPG+$R_{22}$ observamos ainda a ocorrência  de pequenas oscilações na vizinhança do choque, Figura 2d, quando $h =1/3000$. Estas oscilações persistem quando  utilizamos $h=1/4000$, como pode ser confirmado observando a Tabela 3, onde apresentamos os cálculos dos logaritmos dos erros entre as formulações estudadas. Uma alterantiva para amenizar estas oscilações seria utilizar um estabilizador no MEFMQ.

Por fim, verificamos que ao utilizarmos  $\Delta t$ da ordem de $10^{-4}$ as formulações MEFG+$R_{22}$ e SUPG+$R_{22}$ possibilitaram resultados
mais apurados quando comparadas às demais formulações semi-discretas, independente do refinamento do espaço.

\section{Conclusão}
O trabalho forneceu uma comparação entre as formulações semi-discretas para resolver a equação 1D de Burgers, para diferentes condições iniciais e de fronteira. Observamos que o aproximante de Padé $R_{22}$ adicionado as formulações MEFMQ, MEFG e SUPG apresentou soluções aproximadas com convergência mais rápida e maior precisão em comparação ao $R_{11}$. Também, o método $R_{22}$ amenizou consideravelmente as oscilações quando usado nas formulações MEFG e SUPG.

Verificamos no Exemplo 1 que melhorando o refinamento na direção espacial a formulação MEFMQ+$R_{11}$ mostrou-se mais precisa entre as demais formulações, não havendo  a necessidade de um método oneroso para a discretização no tempo. Enquanto que no Exemplo 2 as formulações MEFG+$R_{22}$ e SUPG+$R_{22}$ possibilitaram resultados mais apurados quando comparadas às outras formulações semi-discretas, independente do refinamento do espaço. Salientamos então que, para os exemplos apresentados, predominou a escolha da formulação para a discretização temporal e espacial.

Notamos ainda que os erros obtidos nas formulações semi-discretas, ao refinarmos a malha, diminuiram consideravelmente, apresentando concordância entre as soluções numéricas com a solução exata independente da formulação utilizada.

\begin{abstract}
{\bf Abstract}.
 In this work we compare semi-discrete formulations to obtain
numerical solutions for the 1D Burgers equation. The formulations consist
in the discretization of the time-domain via multi-stage methods of second
and fourth order: $R_{11}$  and $R_{22}$ Padé approximants, and of the spatial-domain
via finite element methods: least-squares (MEFMQ), Galerkin (MEFG) and
Streamline-Upwind Petrov-Galerkin (SUPG). Knowing the analytical solutions
of the 1D Burgues equation, for different initial and boundary conditions,
analyzes were performed for numerical errors from $L_{2}$ and $L_{\infty}$ norm. We
found that the $R_{22}$ Padé approximants, added to the MEFMQ, MEFG, and SUPG
formulations, increased the region of convergence of the numerical
solutions, and showed greater accuracy when compared to the solutions
obtained by the $R_{11}$ Padé approximants. We note that the $R_{22}$ Padé
approximants softened the oscillations of the numerical solutions
associated to the MEFG and SUPG formulations.
\end{abstract}


\begin{thebibliography}{8}

%\bibitem{CLBG} C.S.Q. Caldas, J. Limaco, R.K. Barreto, P. Gamboa, About the
%Benjamin-Bona-Mahony equation in domains with moving boundary,
%{\em TEMA - Tend. Mat. Apl. Comput.}, {\bf 8}, No. 3 (2007), 329--339.
\bibitem{Daryoush e Encyeh}
D. Behmardi,  D.E. Nayeri, Introduction of Fréchet and Gâteaux Derivative. \textit{Appl. Math. Sci.}, \textbf{2} (2008) 975-980.

\bibitem{Brenner e Scott}
C.S. Brenner, R.L. Scott,  ``The Mathematical Theory of Finite Element Methods'', New York, Springer-Verlag, 2008.

\bibitem{Brooks}
A.N. Brooks, T.J.R. Hughes, Streamline upwind/Petrov-Galerkin formulations for convection dominated flows with particular emphasis on the incompressible Navier-Stokes equation,\textit{Comput. Meth. Appl. Mech. Eng.},
\textbf{32} (1982) 199-259.

\bibitem{Ciarlet}
G.P. Ciarlet,  ``The Finite Element Method for Elliptic Problems''. North Holland, SIAM, 1978.

\bibitem{David}
 E.A. David, B. Oscar, Time stepping via one-dimensional Padé approximation, \textit{J. Sci. Comput.}, \textbf{30} (2005), 83-115.

\bibitem{Dogan}
A.A. Dogan, T.J.R. Hughes, A Galerkin element approach to Burgers' equations,\textit{Applied Mathematics
and Computation}, \textbf{154} (2004) 331-346.

\bibitem{DoneaHuerta}
J. Donea, B. Roig, A. Huerta, ``Finite Element Methods for Flow Problems''. John Wiley and Sons, Chichester, 2003.

\bibitem{Donea}
J. Donea, B. Roig,  A. Huerta, Higher-order accurate time-stepping schemes for convection-diffusion problems, \textit{Comput. Meth. Appl. Mech. Engng.}, \textbf{182} (2000) 249-275.

\bibitem{Valdemir1}
V.G. Ferreira, G.A.B. Lima, L. Corrêa, A.C. Cansezano, E.R. Cirilo, P.L. Natti, N.M.L. Romeiro, Avaliação computacional, de esquemas convectivos em problemas de dinâmica dos fluidos. \textit{Semina: Ciências Tecnológicas}, \textbf{2} (2012) 107-116.

\bibitem{Valdemir2}
V.G. Ferreira, R.A.B. de Queiroz, G.A.B. Lima, R.G. Cuenca, C.M. Oishi, J.L.F. Azevedo, S. McKee, A bounded upwinding scheme for computing convection-dominated transport problems, \textbf{57} (2012) 208-224.

\bibitem{Hairer}  E. Hairer, S.P. Norsett, G. Wanner, ``Solving ordinary differential equations I, Non-stiff Problems'', New York, Springer-Verlag, 1987.

\bibitem{Huertabernadino}
A. Huerta, B. Roig, J. Donea, Time-accurate solution of stabilized convection-diffusion-reaction equations: II - accuracy analysis and examples, \textit{Commun. Numer. Meth. Eng.}, \textbf{18} (2002) 575-584.

\bibitem{Jain}
P.C. Jain, R. Shankar, T.V. Singh, Numerical Technique for Solving Convection-Reaction-Diffusion Equation, \textit{Math. Comput. Model}, \textbf{22} (1995) 113-125.

\bibitem{Jiang}
B.N, Jiang, ``The Least-Squares Finite Element Method: Theory and Applications in Computational Fluid Dynamics and Electromagnetics''. Berlin: Springer, 1998.

\bibitem{Kakuda}
K. Kakuda, N. Tosaka,  The generalised boundary element approach to Burgers equation, \textit{Int. J. Numer. Meth. Eng.}, \textbf{29} (1990) 245-261.

\bibitem{Burgers4}
S. Kutluay, A. Esen, I. Dag,  Numerical solutions of the Burgers equation by the least-squares quadratic B-spline finite element method, \textit{J. Comput. Appl. Math.}, \textbf{167} (2004) 21-33.

\bibitem{Cibele}
C.A. Ladeia, ``Formulação semi-discreta aplicada as equações 1D de convecção-difusão-reação e de Burger'', Dissertação de Mestrado, PGMAC/UEL, Londrina, Pr, 2012.

\bibitem{Lambert}
J.D. Lambert, ``Numerical Methods for Ordinary Differential Systems'', New York, Wiley,  1993.

\bibitem{Oden}
J.T. Oden, T. Belytschko, I. Babuska, T.J.R. Hughes, Research directions in computational mechanics, \textit{Comput. Methods Appl. Mech. Eng.}, \textbf{192} (2003) 913-922.

\bibitem{RomeiroActa}
S.R. Pardo, P.L. Natti, N.M.L.  Romeiro, E.R. Cirilo, A transport modeling of the carbon-nitrogen cycle at Igapó I Lake - Londrina, Paraná State, Brazil,
 \textit{Acta Scientiarum, Technology}, \textbf{2} (2012) 217-226.

\bibitem{RomeiroLandau}
N.M.L.  Romeiro, R.S.G Castro, S.M.C. Malta, L. Landau, A linearization technique for multi-species transport problems, \textit{Trans. Porous Med.}, \textbf{70} (2007) 1-10.

\bibitem{Burgers6}
P.L. Sachdev, ``Nonlinear Diffusive Waves'', Cambridge University Press, Cambridge,  1987.

\bibitem{John}
J.C. Strikwerda, ``Finite Difference Schemes and Partial Differential Equations'', SIAM, 2004.

\bibitem{Tian}
Z.F. Tian, P.X. Yu, A High-order exponencial scheme for solving 1D unsteady convection-difusion equations, \textit{Jornal of Computational and Applied Mathematics}, \textbf{235} (2011), 2477-2491.

\bibitem{venutilli}
M. Venutelli, Time-stepping Padé-Petrov-Galerkin models for hydraulic jump simulation. \textit{Math. Comput. Simul.}, \textbf{66} (2004) 585-604.


\end{thebibliography}
\end{document}